
\input amstex.tex
\documentstyle{amsppt}
\magnification=\magstep1
\hsize=12.5cm
\vsize=18cm
\hoffset=1cm
\voffset=2cm
\def\DJ{\leavevmode\setbox0=\hbox{D}\kern0pt\rlap
{\kern.04em\raise.188\ht0\hbox{-}}D}
\footline={\hss{\vbox to 2cm{\vfil\hbox{\rm\folio}}}\hss}
\nopagenumbers
\font\ff=cmr8

\baselineskip=13pt
\def\hf{{\textstyle{1\over2}}}
\def\a{\alpha}
\def\d{{\,\roman d}}
\def\e{\varepsilon}

\def\g{\gamma} 
\def\k{\kappa}

\def\={\;=\;}

\def\zt{\zeta(\hf+it)}

\def\no{\noindent}
\def\R{\Re{\roman e}\,}  
\def\z{\zeta}
\def\no{\noindent}
\def\e{\varepsilon}

\def\no{\noindent}
\def\e{\varepsilon}
\def\l{\lambda}
\def\no{\noindent}
\font\teneufm=eufm10
\font\seveneufm=eufm7
\font\fiveeufm=eufm5
\newfam\eufmfam
\textfont\eufmfam=\teneufm
\scriptfont\eufmfam=\seveneufm
\scriptscriptfont\eufmfam=\fiveeufm
\def\mathfrak#1{{\fam\eufmfam\relax#1}}

\font\tenmsb=msbm10
\font\sevenmsb=msbm7
\font\fivemsb=msbm5
\newfam\msbfam
\textfont\msbfam=\tenmsb
\scriptfont\msbfam=\sevenmsb
\scriptscriptfont\msbfam=\fivemsb
\def\Bbb#1{{\fam\msbfam #1}}

\def \NN {\Bbb N}

\def \RR {\Bbb R}

\def\rightheadline{{\hfil{\ff
Small values  of $|\zt|$}\hfil\tenrm\folio}}

\def\leftheadline{{\tenrm\folio\hfil{\ff
Aleksandar Ivi\'c }\hfil}}
\def\emptyheadline{\hfil}
\headline{\ifnum\pageno=1 \emptyheadline\else
\ifodd\pageno \rightheadline \else \leftheadline\fi\fi}

\topmatter
\title ON SMALL VALUES OF THE RIEMANN ZETA-FUNCTION
ON THE CRITICAL LINE AND GAPS BETWEEN ZEROS\endtitle
\author   Aleksandar Ivi\'c \endauthor
\dedicatory{Lietuvos Matematikos Rinkinys, 42(2002), 31-45}
\enddedicatory
\address{
Aleksandar Ivi\'c, Katedra Matematike RGF-a
Universiteta u Beogradu, \DJ u\v sina 7, 11000 Beograd,
Serbia (Yugoslavia).}
\endaddress
\email {\tt aleks\@ivic.matf.bg.ac.yu,
aivic\@rgf.rgf.bg.ac.yu} \endemail
\abstract
{Small values of $|\zt|$ are investigated, using the value distribution
results of A. Selberg. This gives an asymptotic formula
for $$
\mu\left(\{ 0 < t \le T : |\zt| \le c \}\right).
$$
Some related problems involving gaps between ordinates of zeros of $\z(s)$
are also discussed.}
\endabstract
\endtopmatter

\vglue 1cm

\no
The aim of this note it to discuss the problem of ``small" values of the
Riemann zeta-function $\zeta(s)$ on the critical line
$\Re{\roman e}\, s = \hf$, and some related problems involving the
gaps between the zeros of $\zeta(s)$.  This is
in contrast with the so-called ``large" values of $|\zt|$ (i.e., values
which are $\ge t^\e$), which are extensively discussed in [5]. Since
we have (see [5])
$$
\int_0^T|\zt|^2\d t \sim T\log T,\quad
\int_0^T|\zt|^4\d t \sim {T\over2\pi^2}\log^4 T \quad(T\to\infty),
$$
this means that $|\zt|$ is small ``most of the time". The problem, then,
is to evaluate asymptotically the measure of the subset of $\,[0,T]\,$
where $|\zt|$ is ``small".
\medskip
There are several ways in which one can proceed, and a natural way is
the following one.
Let $c > 0$ be a given constant, let $\mu(\cdot)$ denote measure, and let
$$
A_c(T) := \{ 0 < t \le T : |\zt| \le c \}.
$$
In [8] I raised the question of the asymptotic evaluation of
 $\mu(A_c(T))$. One can tackle this problem by using
the limit law
$$
\lim_{T\to\infty}{1\over T}\,\mu\left(\,\{ 0 < t \le T : |\zt| \le
e^{y\sqrt{{1\over2}\log\log T}}\,\}\right) = {1\over\sqrt{2\pi}}
\int_{-\infty}^y e^{-{1\over2}u^2}\,\d u,\leqno(1)
$$
where $y \in \RR$  is fixed. This result was proved by A. Laurin\v cikas
[11], who used  the fact that
$$
{1\over T}\int_0^T |\zt|^{2k(2\log\log T)^{-1/2}}\d t = e^{{1\over2}k^2}
\left\{1 + O\Bigl((\log\log T)^{-1/4}\Bigr)\right\}\leqno(2)
$$
uniformly for $e^{-\sqrt{\log\log T}} \le k \le k_0$, where $k_0 \in\NN$
is a constant. The proof uses the property that $e^{{1\over2}k^2}$ is
exactly the $2k$--th moment of the distribution function
$$
G(x)  \= \Phi(\log x)  \quad(x > 0),
\quad \Phi(x) := {1\over\sqrt{2\pi}}\int_{-\infty}^x
e^{-{1\over2}u^2}\d u,
$$
so that (2) yields (1) ($G(x) = 0$ for $x \le 0$).
One does not see, however, how one can obtain (1)
from Laurin\v cikas' proof in the form which would not give the result
only as ``lim", but an asymptotic formula with an error term as $T \to
\infty$. This is because the  lognormal   law $G(x)$ is ``bad".
It is known from probability theory that the function
$G(x)$ cannot be defined by its moments $e^{{1\over2}k^2}$. Namely the
moments $e^{{1\over2}k^2}$ are very
rapidly increasing, and from this all ``bad" consequences follow. To obtain
the estimate of the rate of convergence we must consider complex moments,
which one may write as
$$
{1\over
T}\int_0^T|\zt|^{2i\tau(2\log\log
T)^{-1/2}}\d t\;=\;e^{-{1\over2}\tau^2}+S_T(\tau)\qquad(\tau\in\RR),
$$
say. However, the  problem of the estimation of
the function $S_T(\tau)$) seems to be very hard.

\medskip
We shall first show how to  use (1) to obtain a weak asymptotic
formula for $\mu(A_c(T))$.
Let $\e > 0$ be fixed. Note that, for $T \ge T_0(\e,c)$, we trivially have
$$
e^{-\e\sqrt{{1\over2}\log\log T}} < c < e^{\e\sqrt{{1\over2}\log\log T}}.
$$
Therefore, as $T \to \infty$, (1) gives
$$
\eqalign{
\mu(A_c(T)) &\le \mu\bigl(\{
0 < t \le T : |\zt| \le e^{\e\sqrt{{1\over2}\log\log T}}\,\}\bigr)
\cr&
= {1\over\sqrt{2\pi}}
\int\limits_{-\infty}^\e e^{-{1\over2}u^2}\d u\cdot
T + o(T),\cr}
$$
and
$$
\eqalign{
\mu(A_c(T)) &\ge \mu\bigl(\{
0 < t \le T : |\zt| \ge e^{-\e\sqrt{{1\over2}\log\log T}}\,\}\bigr)\cr&
= {1\over\sqrt{2\pi}}
\int\limits_{-\infty}^{-\e} e^{-{1\over2}u^2} \d u\cdot
T + o(T).\cr}
$$
But as
$$
\int_{-\infty}^0 e^{-{1\over2}u^2}\d u = \sqrt{2}\,\int_0^\infty
e^{-x^2}\d x = \sqrt{\pi\over2},
$$
it follows that
$$
\mu(A_c(T)) = {T\over2} + O(\e T) + o(T),
$$
hence letting $\e \to 0$ we obtain
$$
\mu(A_c(T)) = {T\over2} + o(T)\qquad(T\to\infty).\leqno(3)
$$

\medskip

Let now $0 < c_1 < c_2$. Since
$$
[c_1, c_2] = (0,c_2] \;\backslash\; (0,c_1],
$$
it follows from (3), as $T \to \infty$, that
$$
\mu\left(\{ 0 < t \le T : c_1 \le |\zt| \le c_2 \}\right)
= \mu(A_{c_2}(T)) - \mu(A_{c_1}(T)) + o(T) = o(T). \leqno(4)
$$

\medskip
It turns out that for the above problems one can
use Theorem 2 of A. Selberg's paper [13], which is an asymptotic
formula with an error term. Selberg obtained sharper results than
Laurin\v cikas' before  Laurin\v cikas did, but he published his paper later.
Actually Selberg's paper contains no proofs, but it is hinted at the
end that proofs  will appear. Also there exists the recent work of
D.A. Hejhal [4], which is built on the methods of [13] and complements it.
In fact (2.6) of Theorem 2 on p. 374 of Selberg's paper can be
specialized to yield a result sharper than (1), namely
$$\eqalign{
&\mu\left(\,\Big\{ 0 < t \le T : |\zt| \le
e^{y\sqrt{{1\over2}\log\log T}}\,\Big\}\right)\cr&
= \Phi(y)T + O\left(T{(\log\log \log T)^2\over\sqrt{\log\log T}}\right),}
\leqno(5)$$
where as before, for $x \in \RR$,
$$
\Phi(x) = {1\over\sqrt{2\pi}}\int_{-\infty}^x
e^{-{1\over2}u^2}\d u
$$
is the probability integral. Now Hejhal kindly confirmed, by
going through Selberg's unpublished proof, that
formula (5) holds
{\bf uniformly} in $y$. Therefore choosing
$$
y \;=\; {\log c\over\sqrt{{1\over2}\log\log T}}
$$
for a given constant $c > 0$, and using the fact that, for $|y| \le 1$,
$$
\Phi(y) = {1\over\sqrt{2\pi}}\int_{-\infty}^0
e^{-{1\over2}u^2}\d u + O\left(\int_0^{|y|}
e^{-{1\over2}u^2}\d u\right) = {1\over2} + O(|y|),
$$
we obtain from (5)

\bigskip
{\bf THEOREM 1.} {\it We have}
$$
\mu(A_c(T)) = {T\over2} +
O\left(T{(\log\log \log T)^2\over\sqrt{\log\log T}}\right),\leqno(6)
$$
{\it where as before}
$$
A_c(T) := \{ 0 < t \le T : |\zt| \le c \}.
$$
{\it We also have, for given constants $0 < c_1 < c_2$,}
$$
\mu\left(\Big\{ 0 < t \le T : c_1 \le |\zt| \le c_2 \Big\}\right)
= O\left(T{(\log\log \log T)^2\over\sqrt{\log\log T}}\right).\leqno(7)
$$

\bigskip\no
Of course, (7) follows easily from (4) and (6). We note that the
formulas (6) and (7), which improve (3) and (4), give a
satisfactory solution to the problem of
the distribution of ``small" values of $|\zt|$. The factor
$(\log\log \log T)^2$, which appears in (5)--(7), is probably extraneous,
but will be very likely difficult to get rid of.

Another way to see how (6) and (7) follow is to apply a result
contained in D.A. Hejhal's work [4], where he successfully deals with
zeros of linear combinations of $L$-functions belonging to Selberg's
class [13]. In particular, his equation (4.21), specialized to
$\zeta(s)$, says that
$$
\mu\left(\{\,T \le t \le 2T \,:\, e^a \le |\zt| \le e^b\,\}\right)
= T\int_{a/\sqrt{\pi\psi}}^{b/\sqrt{\pi\psi}}\,e^{-\pi v^2}\d v
+ O\left({T\log^2\psi\over\sqrt{\psi}}\right)\leqno(8)
$$
uniformly in $a,b \in \RR$, where
$$
\psi \= \log\log T + O(\log\log\log T).
$$
Therefore the specialization
$a = -\infty,\, b = \log c$ yields (6), while $a = \log c_1,\,
b =\log c_2\,(0 < c_1 < c_2)$ yields (7).
\bigskip

\medskip

We shall  consider now a problem related to the above one. Let henceforth
$0 < \g_1 \le \g_2 \le \ldots\,$ denote positive ordinates of complex
zeros of $\z(s)$; it is known that $\g_1 = 14.13\ldots\,$, and all
known ($> 10^9$) zeros are simple and lie on the critical line
$\R s = \hf$.  We
define $\g_-(t) = \g_n$ if $\g_n \le t < \g_{n+1},\, \g_+(t) = \g_{n+1}$
if $\g_n < t \le \g_{n+1},\, \g_-(t) = \g_+(t) = \g_n$ if $t = \g_n$,
$$
{\Cal A}(T) \= \{\,0 < t \le T\,:\,|\zt| \le \g_+(t) - \g_-(t)\,\},
$$
$$
{\Cal B}(T) = [0,\,T] \;\backslash\; {\Cal A}(T)
\=  \{\,0 < t \le T\,:\,|\zt| > \g_+(t) - \g_-(t)\,\}.
$$
Natural problems are  to evaluate asymptotically
$\mu({\Cal A}(T))$ and $\mu({\Cal B}(T))$. We shall prove
the following

\bigskip
{\bf THEOREM 2.}  {\it We have}
$$
\mu({\Cal B}(T)) \= T +
O\left(T{\log\log \log T\over\sqrt{\log\log T}}\right).\leqno(9)
$$

\bigskip\no
{\bf Proof.}  We shall first employ a method based on the
value distribution result (8). This leads to (9), but with
$(\log\log\log T)^2$ in place of $\log\log \log T$. Then we shall
present another approach, which yields the slightly sharper
result of Theorem 2. Let
$$
{\Cal C}_1(T) := \left\{\,0 < t \le T\,:\, \g_+(t) - \g_-(t) <
{(\log\log T)^6\over\log T}\,\right\},
$$
$$
{\Cal C}_2(T) := \{\,0 < t \le T\,:\,|\zt| > \exp(-(\log\log T)^{3/4})\,\},
$$
and let $\bar S$ denote the complement of $S$ in $[0,\,T]$.
From (8) with $a = -\infty,\, b = -(\log\log T)^{3/4}$ we obtain
$$
\mu({\bar{\Cal C}}_2(T)) \ll T\int_{{1\over2}(\log\log T)^{1/4}}^\infty
e^{-\pi v^2}\d v + T{(\log\log \log T)^2\over\sqrt{\log\log T}}
\ll T{(\log\log \log T)^2\over\sqrt{\log\log T}},
$$
hence
$$
\mu({\Cal C}_2(T)) = T + O\left(T{(\log\log \log T)^2\over\sqrt{\log\log T}}
\right).\leqno(10)
$$
On the other hand
$$
\mu({\Cal C}_2(T)) = \mu({\bar{\Cal C}}_1(T)\cap{\Cal C}_2(T)) +
\mu({\Cal C}_1(T)\cap{\Cal C}_2(T)).\leqno(11)
$$
However we have
$$
\mu({\bar{\Cal C}}_1(T)\cap{\Cal C}_2(T)) \le \mu({\bar{\Cal C}}_1(T))
\ll T{(\log\log \log T)^2\over\sqrt{\log\log T}}.\leqno(12)
$$
The second bound in (12) is a consequence of a bound which follows
from the following Lemma (weaker results are given in A. Fujii [1], [2]
and (without proof)  in E.C. Titchmarsh [14, p. 246]).

\bigskip
{\bf Lemma.} {\it Let $0 < \g_1 \le \g_2 \le \cdots $ denote
imaginary parts of complex zeros of $\z(s)$, and let $\l \ge 2$. Then
there exists a constant $C > 0$ such that uniformly}
$$
\sum_{T<\g_n\le T+H,\g_{n+1}-\g_n\ge\l/\log T}1
\ll (N(T+H)-N(T))\exp\left(-C\l\right) + 1,\leqno(13)
$$
{\it where $N(T)$ is the number of zeros of $\z(s)$ with imaginary
parts in $(0,T]$, and $T^a < H \le T,\, a > \hf$}.

\bigskip\no
{\bf Proof.} The basic result is the asymptotic formula [15, Theorem 4]
of K.-M. Tsang. This says that, for $T^a < H \le T,\, a > \hf$, $0 < h < 1$
and any $k \in\NN$, we have uniformly
$$
\eqalign{
\int_T^{T+H}(S(t+h) - S(t))^{2k}\d t &= {H(2k)!\over(2\pi^2)^kk!}
\log^k(2+h\log T) \cr&
+ O\left\{H(ck)^k\left(k^k + \log^{k-{1\over2}}(2+h\log T)\right)\right\},
\cr}\leqno(14)
$$
where $c > 0$ is a constant, and as usual $S(T) = {1\over\pi}\arg\zeta(
\hf + iT)$.
Thus $S(T) = O(\log T)$ (see [5] or [14])  and the Riemann--von Mangoldt
formula is
$$
N(T) \= {T\over2\pi}\log\left({T\over2\pi}\right) - {T\over2\pi}
+ S(T) + {7\over8} + O\left({1\over T}\right).
$$
This gives $\g_{n+1} - \g_n \ll 1$, and also
$$
S(t+h) - S(t) = N(t+h) - N(t) - {h\over2\pi}\log t + O\left(
{h^2 +1\over t}\right).\leqno(15)
$$
If
$$
\g_n < t < \hf(\g_n + \g_{n+1}),\, \g_{n+1} - \g_n \ge {\l\over\log T},
\, T \le t \le T+H,\, h = {\l\over2\log T},\leqno(16)
$$
then
$N(t+h) - N(t) = 0$, and $h \ll 1$ will hold in view of
$\g_{n+1} - \g_n \ll 1$. For $t$ satisfying (16) we have
$$
|S(t+h) - S(t)| \;\ge\; {h\over4\pi}\log t \;\ge\; {\l\over8\pi},
$$
and (14) will in fact hold for $0 < h \ll 1$. We obtain from (14)
$$
\sum_{T<\g_n<\g_{n+1}\le T+H,\g_{n+1} - \g_n \ge \l/\log T}
\left({\l\over8\pi}\right)^{2k}(\g_{n+1} - \g_n)
\ll H(Ak(k+\log\l))^k
$$
with suitable $A > 0$, which implies that ($B = (8\pi)^2A$)
$$
\sum_{T<\g_n\le T+H,\g_{n+1} - \g_n \ge \l/\log T}1
\ll (N(T+H) - N(T))\left(Bk{(k+\log\l)\over\l^2}\right)^k + 1.\leqno(17)
$$
We take
$$
k \;=\; \left[{\l\over2\sqrt{B}}\right],
$$
and (13) follows from (17)
for $\l\ge\l_0\,(\ge 2)\,$, while for $\l < \l_0$ the bound in (13)
is trivial.

\medskip

To obtain (12) write
$$
{\bar{\Cal C}}_1(T) = \bigcup_{k=1}^\infty D_k(T),
D_k(T) := \left\{0 < t \le T: V_k(T) \le \g_+(t) - \g_-(t)
< 2V_k(T)\right\},
$$
$$
V_k(T) \;:=\; {2^{k-1}(\log\log T)^6\over \log T}.
$$
Hence with $\lambda = \lambda(k,T) = 2^{k-1}(\log\log T)^6$ we have,
on using (13),
$$\eqalign{
\mu(D_k(T)) &\le 2V_k(T)\sum_{\g_n \le T, \g_{n+1}
- \g_n \ge \lambda/\log T}1\cr&
\ll T\exp(-2^{k}(\log\log T)^2),\cr}
$$
which gives
$$
\mu({\bar{\Cal C}}_1(T)) \ll T\sum_{k=1}^\infty
\exp(-2^{k}(\log\log T)^2) \ll
T{(\log\log \log T)^2\over\sqrt{\log\log T}},
$$
as asserted.

We therefore have from (10)--(12)
$$
\mu({\Cal C}_1(T)\cap{\Cal C}_2(T)) = T + O\left(
T{(\log\log \log T)^2\over\sqrt{\log\log T}}\right),
$$
and (9) with the error  term
$O\left(T{(\log\log \log T)^2\over\sqrt{\log\log T}}\right)$
follows from
$$
T \ge \mu({\Cal B}(T)) \ge \mu({\Cal C}_1(T)\cap{\Cal C}_2(T)),
$$
since for $t \in {\Cal C}_1(T)\cap{\Cal C}_2(T)$ we have
$$
|\zt| > e^{-(\log\log T)^{3/4}} > {(\log\log T)^6\over \log T}
> \g_+(t) - \g_-(t).
$$

\bigskip
To obtain Theorem 2 in the sharper form given by (9), we use
a result of A. Perelli and the author [10] (see also [6, Theorem 6.2]) which
says that, if $\psi(T)$ is an arbitrary positive function tending
to infinity with $T$, then for
$$
0 \le \l \le \left(\psi(T)\log\log T\right)^{-1/2}
$$
we have
$$
\int_0^T|\zt|^\l \d t = T + o(T)\qquad(T \to \infty).\leqno(18)
$$
For our purposes we need (18) with an $O$-term for the error instead
of the $o$-term. This is given by
$$
\int_0^T|\zt|^\l \d t = T + O\left({T\over\sqrt{\psi(T)}}\right)
+ O(\log T).\leqno(19)
$$
To obtain (19) in place of (18) one has first to note that [6, Lemma 6.7]
actually gives
$$
\int_0^T|\zt|^\l \d t \ge T + O(\log T)
$$
for any $\l \ge 0$. For the corresponding upper bound it suffices to
note that, in the proof of [6, (6.41)] we obtain, for
$m = [(\log\log T)^{1/2}],\, C_1 > 0,\, C_2 > 0$,
$$\eqalign{
(C_1)^{{1\over2}m\l}(\log T)^{\l\over2m}\ &\le \exp\left(
\left(C_2\psi(T)\log\log T\right)^{-1/2}(\log\log T)^{1/2}\right)\cr&
=  1 + O\left({1\over\sqrt{\psi(T)}}\right),\cr}
$$
which gives then  (19), as asserted. Let henceforth
$$
\l \;:=\; {1\over\sqrt{\psi(T)\log\log T}},
\quad
\psi(T) \;:=\; {\log\log T\over9(\log\log\log T)^2}.
$$
On one hand, we have  (19), while on the other hand we may write
$$
\int_0^T|\zt|^\l \d t = \int_{{\Cal A}(T)}|\zt|^\l \d t
+ \int_{{\Cal B}(T)}|\zt|^\l \d t = I_1(T) + I_2(T),\leqno(20)
$$
say. For $I_2(T)$ we use the Cauchy-Schwarz inequality and (19) with
$2\l$ replacing $\l$ to obtain that
$$\eqalign{
I_2(T) &\le (\mu({\Cal B}(T))^{1/2}\left(\int_0^T|\zt|^{2\l}
\d t\right)^{1/2}\cr&
= (\mu({\Cal B}(T))^{1/2}\left(T + O\left({T\over\sqrt{\psi(T)}}\right)
\right)^{1/2}.\cr}\leqno(21)
$$
We have
$$
I_1(T) \;\le\;\sum_{\g_n\le T}(\g_{n+1}-\g_n)^{\l+1}
\ll T\log^{-\l}T = T(\log\log T)^{-3}.\leqno(22)
$$
Here we used the bound
$$
\sum_{\g_n\le T}(\g_{n+1}-\g_n)^{\a} \ll T(\log T)^{1-\a}\qquad
(1 \le \a \le \a_0).\leqno(23)
$$
To obtain (23) we estimate trivially the contribution of
$\g_n$ for which  $\g_{n+1}-\g_n \le 2/\log T$. The remaining sum is
split into subsums where
$$
{2^k\over\log T} < \g_{n+1}-\g_n \le {2^{k+1}\over\log T}
\qquad(k = 1,2,\ldots\,),
$$
each of which is estimated by (13), which yields (23).

Therefore we obtain  from (15)--(22) a lower bound for $\mu({\Cal B}(T))$
of the form given by (9), and trivially
$\mu({\Cal B}(T)) \le T$. This establishes (9).

\bigskip
It is very likely that
preceding results hold if the $\g_n$'s are the ordinates of
zeros on the critical line (assuming that the Riemann hypothesis is
not true, and $\z(s)$ has zeros lying off the critical line), but in that
case the problems are more difficult. Connected with this is
a problem which I posed during the Conference on
Elementary and Analytic Number Theory, held in Oberwolfach, March 1994
(see also [7]). This is also related to small values of $|\zt|$.
Let ${\bar \g}_n$ denote the $n$-th positive zero of $\zeta({1\over2} + it)
 = 0$, where possible
multiple zeros are counted with their respective multiplicities. Let
$$
N_0(T) = \sum_{{\bar \g}_n\le T}1,\qquad B(T) \;:=\; N_0(T) - A(T),
$$
$$
A(T) \;:=\; \sum_{{\bar \g}_n\le T,\max\limits_{{\bar \g}_n\le t\le
{\bar \g}_{n+1}}
|\zeta({1\over2}+it)|\le{\bar \g}_{n+1}-{\bar \g}_n}1.
$$
The problem is to compare (unconditionally, or under the Riemann
hypothesis)
$A(T)$ and $B(T)$ to $N_0(T)$ (we know that
$T\log T \ll N_0(T) \ll T\log T$).
I expect that $B(T) \sim N_0(T)$ (or equivalently
$A(T) = o(N_0(T))$) as $T\to\infty$, that is, on the
average the maximum between two consecutive zeros on the critical
line should be larger than the gap between these zeros. M. Jutila and
the author [9] proved that the number of ${\bar \g}_n$'s not exceeding $T$
for which ${\bar \g}_{n+1} - {\bar \g}_n \ge V \, (> 0)$ is uniformly
$$
\ll \quad \min(TV^{-2}\log T,\,TV^{-3}\log^5T),
$$
but unfortunately this bound is not well
suited in dealing with the ``small gaps".

\bigskip Returning to Theorem 2,
note that ${\Cal A}(T)$ contains intervals $\,[\g_n,\,\g_{n+1}]\,$ with
$\g_n \le T$ (with the possible exception of one interval), such that
$$
\max_{\g_n\le t\le\g_{n+1}}|\zt| \;\le\; \g_{n+1} - \g_n.
$$
Then the method of proof of Theorem 2 shows that
$$
{\sum_{\g_n\le T}}^* (\g_{n+1} - \g_n) = T + O\left(
T{\log\log \log T\over\sqrt{\log\log T}}\right),\leqno(24)
$$
where ${}^*$ denotes summation with the conditions
$$
\g_{n+1} - \g_n < {(\log\log T)^6\over\log T},\quad
\max_{\g_n\le t\le\g_{n+1}}|\zt| > \g_{n+1} - \g_n.
$$
Now we assume the Riemann hypothesis and
apply the Cauchy-Schwarz inequality to the left-hand side
of (24). Then by (23) with $\a = 2$ we obtain
$$
B(T) \;\gg\; N_0(T),\leqno(25)
$$
which favours the conjecture that $B(T) \sim N_0(T)$ as $T\to\infty$.
Actually the constant in (25) may be explicitly calculated if we  use
a bound of A. Fujii [3], namely
$$
\sum_{{\bar \g}_n\le T}({\bar \g}_{n+1} - {\bar \g}_n)^2 \;\le\;
9\cdot{2\pi T\over\log{T\over2\pi}}\qquad(T > T_0).\leqno(26)
$$
This leads, under RH, to the inequality
$$
B(T) \;\ge\; (1 + o(1)){T\over18\pi}\log\left({T\over2\pi}\right) \;=\;
({1\over9} + o(1))N_0(T)\qquad(T\to\infty).
$$
If, in addition to the RH, one assumes the Gaussian
Unitary Ensemble Hypothesis, then one can improve the bound in (26)
and obtain in fact an asymptotic formula for the sum on the left-hand side of
(26). For the details the reader is referred to [5].

\medskip
Note that $A(T)$ trivially counts the ${\bar \g}_n$'s for which
${\bar \g}_n = {\bar \g}_{n+1}$, that is, multiple zeros on the critical
line. Hence the conjecture $B(T) \sim N_0(T)$ is stronger than
the conjecture that almost all zeros on the critical line are simple
(which seems to be independent of the RH). In connection with this
it is perhaps natural to consider also
$$
D(T) \;:=\; \sum_{{\bar \g}_n< {\bar \g}_{n+1}\le T,
\max\limits_{{\bar \g}_n\le t\le {\bar \g}_{n+1}}
|\zeta({1\over2}+it)|\le{\bar \g}_{n+1}-{\bar \g}_n}1,
$$
and try to show that
$$
D(T) \;=\; o(N_0(T))\qquad(T\to\infty),\leqno(27)
$$
which is implied by $A(T) = o(N_0(T))$.
One way to deal with this problem is to note that from Theorem 2 we have,
unconditionally,
$$
\mu({\Cal A}(T)) = \sum_{\g_n\le T,\max\limits_{\g_n\le t\le
{ \g}_{n+1}}|\zeta({1\over2}+it)|\le{\g}_{n+1}-{\g}_n}
(\g_{n+1}-\g_n) \ll T{\log\log \log T\over\sqrt{\log\log T}}.
\leqno(28)
$$
On the other hand, for any $\k > 0$,
$$\eqalign{
\mu({\Cal A}(T)) &\ge
\sum_{\g_n\le T,\max\limits_{\g_n\le t\le
{ \g}_{n+1}}|\zeta({1\over2}+it)|\le{\g}_{n+1}-{\g}_n, \g_{n+1}-\g_n
>\k/\log T}( \g_{n+1}-\g_n) \cr&
\ge {\k\over\log T}
\sum_{\g_n\le T,\max\limits_{\g_n\le t\le
{ \g}_{n+1}}|\zeta({1\over2}+it)|\le{\g}_{n+1}-{\g}_n, \g_{n+1}-\g_n
>\k/\log T}1.\cr}
$$
Therefore from (28) we obtain, as $T\to\infty$,
$$
\sum_{\g_n\le T,\max\limits_{\g_n\le t\le
{ \g}_{n+1}}|\zeta({1\over2}+it)|\le{\g}_{n+1}-{\g}_n, \g_{n+1}-\g_n
>\k/\log T}1
\ll {\log\log \log T\over\k\sqrt{\log\log T}}N(T) = o(N(T))
\leqno(29)
$$
provided that
$$
{\log\log \log T\over\k\sqrt{\log\log T}} \;=\;o(1)
\qquad(T\to\infty).\leqno(30)
$$
Therefore if we can show that
$$
\sum_{\g_n<\g_{n+1}\le T, \g_{n+1}-\g_n \le\k/\log T}1 = o(N(T))
\qquad(\k = o(1),\;T\to\infty)\leqno(31)
$$
for $\k$ satisfying (30), then from (29) and (31) we obtain,
assuming RH, the conjectural relation (27).

However,  (31) follows from what is known as
{\it the essential simplicity hypothesis} of zeta zeros. This says that
$$
\sum_{0<\g,\g'\le T,0<\g-\g'\le2\pi\a/\log(T/2\pi)}1 = o(N(T)) \leqno(32)
$$
for $\a = o(1),\,T\to\infty$ together with the relation
$$
\sum_{0<\g\le T}m(\g) = (1 + o(1))N(T) \qquad(T\to\infty), \leqno(33)
$$
where $\g$ and $\g'$ denote ordinates of zeta zeros, and
$m(\g)$ denotes the multiplicity of the zeta zero $\hf+i\g$ (assuming RH),
which is already counted in the above sum with its muliplicity.
Thus (32) says that pairs of different zeros with small gaps are rare,
while (33) asserts that almost all zeros are simple. In particular, the
essential simplicity hypothesis implies not only (27) but
the stronger $A(T) = o(N_0(T))$ as well (under the RH).
A discussion of the essential simplicity
hypothesis is given by J. Mueller [12]. It is shown there
that this hypothesis is, under the RH, equivalent to two other
hypotheses involving certain integrals. The relation (31) follows as
the limiting case of the Gaussian Unitary Ensemble hypothesis, and
it follows also from the limiting case of Montgomery's
pair correlation conjecture that for fixed $\a > 0$ and $T\to\infty$
$$
\sum_{0<\g,\g'\le T,0<\g-\g'\le2\pi\a/\log(T/2\pi)}1
= \left\{\int_0^\a\left(1-\left({\sin\pi t\over\pi t}\right)^2\right)\d t
+ o(1)\right\}N(T).
$$
Thus proving $A(T) = o(N_0(T))$ (or the weaker (27)) assuming only
the RH seems to be difficult, while an unconditional proof is certainly
out of reach at present.

\bigskip
{\bf Acknowledgement.} I wish to thank A. Fujii,
D.R. Heath-Brown, D.A. Hejhal
and A. Laurin\v cikas for valuable remarks.
\bigskip
\topglue2cm
\Refs

\bigskip

\item{[1]} A. Fujii, {\it On the distribution of zeros of
the Riemann zeta function in short intervals},
Bull. Amer. Math. Soc. {\bf81}(1975), 139-142.

\item{[2]} A. Fujii, {\it On the difference between $r$ consecutive
ordinates of the zeros of the Riemann zeta function}, Proc. Japan
Acad. {\bf51}(1975), 741-743.

\item{[3]} A. Fujii, {\it On the gaps between the consecutive
zeros of the Riemann zeta function}, Proc. Japan Acad. {\bf 66} Ser. A
(1990), 97-100.

\item{[4]} D.A. Hejhal, {\it On a result of Selberg concerning
zeros of linear combinations of L-functions}, International
Math. Res. Notices 2000, No. {\bf11}(2000), 551-557.

\item{[5]} A. Ivi\'c, {\it The Riemann zeta-function}, John Wiley
\& Sons, New York, 1985.

\item{[6]} A. Ivi\'c, {\it Mean values of the Riemann zeta-function},
TATA Institute of Fundamental Research LNs {\bf 82}, Springer,
Berlin etc., 1991.

\item{[7]} A. Ivi\'c, {\it On sums of gaps between the zeros of
$\zeta(s)$ on the critical line}, Univ. Beograd. Publ. Elektrotehn.
Fak. Ser. Mat. {\bf6}(1995), 55-62.

\item{[8]} A. Ivi\'c, {\it Some problems on mean values of the Riemann
zeta-function}, Journal de Th\'eorie des Nombres Bordeaux {\bf 8}(1996),
101-122.

\item{[9]} A. Ivi\'c and M. Jutila, {\it Gaps between consecutive zeros of
the Riemann zeta-function}, Monatshefte Math.  {\bf105} (1988), 59-73.

\item{[10]} A. Ivi\'c and A. Perelli, {\it Mean values of  certain
zeta-functions on the critical line}, Litovskij Mat. Sbornik {\bf29}(1989),
701-714.

\item{[11]} A. Laurin\v cikas, {\it Limit theorem for the
Riemann zeta-function on the critical line II (in Russian)}, Lietuvos
Mat. Rinkinys {\bf27}(1987), 489-500.

\item{[12]} J. Mueller, {\it Arithmetical equivalent of essential
simplicity of zeta zeros}, Trans. Amer. Math. Soc. {\bf 275}(1983), 175-183.

\item{[13]} A. Selberg, {\it Old and new conjectures and results about a
class of Dirichlet series}, in ``Proceedings of the Amalfi Conference
on Analytic Number Theory" (eds. E. Bombieri et al.), Universit\`a di
Salerno, Salerno, 1992,  367-385; also {\it Collected Papers}
(Vol. II), Springer Verlag, Berlin etc., 1991, 47-63.

\item{[14]} E.C. Titchmarsh, {\it The theory of the Riemann
zeta-function}, 2nd edition, Oxford University Press, Oxford, 1986.

\item{[15]} K.-M. Tsang, {\it Some $\Omega$-theorems for the Riemann
zeta-function}, Acta Arith. {\bf46}(1986), 369-395.

\bigskip

Aleksandar Ivi\'c

Katedra Matematike RGF-a

Universitet u Beogradu, \DJ u\v sina 7

11000 Beograd, Serbia and Montenegro

e-mail: aivic\@rgf.bg.ac.yu, aivic\@matf.bg.ac.yu

\endRefs

\bye